\documentclass[conference]{IEEEtran}
\usepackage[utf8]{inputenc}

\usepackage{amsmath,amsfonts,amsthm,amsopn,amssymb}
\usepackage{graphicx}
\usepackage{subfigure}

\usepackage{twoopt}

\usepackage{mathtools}
\mathtoolsset{showonlyrefs}

\usepackage{epsfig,verbatim}
\usepackage{mathrsfs}

\providecommand{\norm}[1]{\lVert#1\rVert}

\newcommand{\R}{\mathbb{R}}
\newcommand{\C}{\mathbb{C}}
\newcommand{\D}{\mathbb{D}}

\newcommand{\Z}{\mathbb{Z}}
\newcommand{\N}{\mathbb{N}}

\newcommand{\Lt}[1][d]{L^2(\R^{#1})}
\newcommand{\G}{\mathcal{G}}

\newcommand{\indicator}{\raisebox{2pt}{$\chi$}}
\renewcommand{\l}{\lambda}
\renewcommand{\L}{\Lambda}

\newcommand{\cond}{\textnormal{cond}}

\theoremstyle{plain}
\newtheorem{theorem}{Theorem}[section]
\theoremstyle{plain}
\newtheorem{conjecture}[theorem]{Conjecture}

\allowdisplaybreaks
\newcommandtwoopt{\xarrow}[2][0.5cm][0]{\mathrel{\rotatebox[origin=c]{#2}{$\xrightarrow{\rule{#1}{0pt}}$}}}

\begin{document}

\title{The Strohmer and Beaver Conjecture for Gaussian Gabor Systems\\--\\A Deep Mathematical Problem (?)}

\author{
	\IEEEauthorblockN{Markus Faulhuber}
	\IEEEauthorblockA{NuHAG, Faculty of Mathematics, University of Vienna\\Oskar-Morgenstern-Platz 1, 1090 Vienna, Austria}
}

\maketitle

\begin{abstract}
	In this article we are going to discuss the conjecture of Strohmer and Beaver for Gaussian Gabor systems. It asks for an optimal sampling pattern in the time-frequency plane, where optimality is measured in terms of the condition number of the frame operator. From a heuristic point of view, it seems obvious that a hexagonal (sometimes called triangular) lattice should yield the solution. The conjecture is now open for 16 years and only recently partial progress has been made. One point this article aims to make, is to show up parallels to a long standing, open problem from geometric function theory, Landau's problem posed in 1929, suggesting that the conjecture of Strohmer and Beaver is a very deep mathematical problem.
\end{abstract}

\IEEEpeerreviewmaketitle

\section{Introduction}
Gaussian Gabor systems have already been studied in the setting of quantum mechanics by John von Neumann in 1932 \cite{Neu32}, but their name originates from the 1946 paper of physics Nobel laureate Dennis Gabor \cite{Gab46}. Gabor was looking for a 2-dimensional representation of a one-dimensional signal (function), providing information on the signal's time and frequency content at the same time. This idea leads to a representation of a signal in a mixed time-frequency domain.

The aim of Gabor analysis is to expand a signal from the Hilbert space $\Lt[]$ into a generalized Fourier series, similar to the Fourier series expansion for functions in $L^2(\mathbb{T})$. Also, the coefficients in the Fourier series should already provide accurate information about the joint time-frequency content of the signal.

Originating from this idea, a whole new field of (harmonic) analysis has been developed -- time-frequency analysis. This field is a very active field of research in mathematics as well as in engineering and has already come up with deep mathematical problems, such as Feichtinger's conjecture. It was mentioned in print for the first time in 2005 in \cite{Cas05} and turned out to be equivalent to the Kadison-Singer conjecture or paving conjecture \cite{KadSin59}, which dates back to 1959. A (positive) solution was finally given in 2015 by Marcus, Spielman and Srivastava \cite{Mar15}, turning the conjectures into theorems.

The aim of this article is to discuss in some detail another conjecture from the field of time-frequency analysis, namely the conjecture of Strohmer and Beaver on optimal Gaussian Gabor frames \cite{StrBea03}, which appeared in print in 2003, and to briefly describe (loose or deep) connections to a problem from geometric function theory, posed by Landau in 1929.

\section{Gabor Systems and Gabor Frames}
Before we can properly define a Gabor system, we need to fix the notation and define some auxiliary tools first. We stick close to the textbook of Gröchenig \cite{Gro01} with our notation. We start with defining the inner product for the Hilbert space $\Lt[]$. The inner product of two functions $f,g \in \Lt[]$ is given by
\begin{equation}
	\langle f, g \rangle = \int_\R f \, \overline{g} \, d \mu,
\end{equation}
where $\overline{g}$ is the complex conjugate of $g$ and $d \mu$ denotes the Lebesgue measure on $\R$. The canonical norm induced by the inner product is
\begin{equation}
	\norm{f}_2^2 = \langle f, f \rangle.
\end{equation}
In the sequel we will usually abuse notation and write out the formulas for point-wise defined functions. The Fourier transform of a function is given by
\begin{equation}
	\widehat{f} (\omega) = \int_\R f(t) e^{-2 \pi i \omega t} \, dt.
\end{equation}
The Fourier transform is unitary on $\Lt[]$, i.e., $\norm{f}_2 = \norm{\widehat{f}}_2$. The fundamental operators in time-frequency analysis are the translation (or time-shift) operator $T_x$ and the modulation (or frequency-shift) operator $M_\omega$. They act on functions by the rules
\begin{equation}
	T_x f(t) = f(t-x) \qquad \textnormal{ and } \qquad M_\omega f(t) = f(t) e^{2 \pi i \omega t}.
\end{equation}
The composition is called a time-frequency shift and denoted by
\begin{equation}
	\pi(\l) = M_\omega T_x, \qquad \l = (x, \omega) \in \R^2.
\end{equation}
In general, time-frequency shifts do not commute, as already the translation and modulation operator do not commute. However, they fulfill the following commutation relation;
\begin{equation}\label{eq_commrel}
	T_x M_\omega = e^{-2 \pi i x \omega} M_\omega T_x.
\end{equation}

For a so-called window $g \in \Lt[]$ and an index set $\L \subset \R^2$, the collection of time-frequency shifted versions of $g$ with respect to $\L$ is called a Gabor system;
\begin{equation}
	\G(g, \L) = \{ \pi(\l) g \mid \l \in \L \}.
\end{equation}
The index set $\L$ is usually a discrete, relatively separated subset of the time-frequency plane $\R^2$. In the so-called regular case, the index set is assumed to be a lattice, which means that it possesses a group structure. In this case, it can be represented by an invertible matrix $M$ and the columns of the matrix serve as a basis for the lattice;
\begin{equation}
	\L = M \Z^2.
\end{equation}

Now, the idea is to expand a function $f \in \Lt[]$ into a series of the form
\begin{equation}\label{eq_series}
	f(t) = \sum_{\l \in \L} c_\l \, \pi(\l) g(t) = \sum_{(x, \omega) \in \L} c_{(x,\omega)} g(t-x) e^{2 \pi i \omega t}.
\end{equation}
A first example of a Gabor system, which actually constitutes an orthonormal basis, is given by
\begin{equation}
	\G(\indicator_{[0,1)}, \Z^2) = \{ \indicator_{[0,1)}(t-k) e^{2 \pi i l t} \mid (k,l) \in \Z^2 \}.
\end{equation}
This system consists of integer shifted copies of the Fourier basis on the torus $\mathbb{T} = \R \slash \Z \cong [0,1)$ along the real line.

However, for various reasons, Gaussian windows are often the preferred choice. One of the many preferable properties of Gaussians is that they uniquely minimize the uncertainty principle. The standard Gaussian window is given by
\begin{equation}
	g_0(t) = 2^{1/4} e^{-\pi t^2},
\end{equation}
where the factor in front is to normalize the Gaussian, i.e., $\norm{g_0}_2 = 1$. Also, throughout this work we will always assume that the window is normalized.

As stated earlier, the system $\G(g_0, \Z^2)$ was already studied by von Neumann \cite{Neu32} and later again by Gabor \cite{Gab46}. As we know now, this system is complete, but there is no stable way to expand a signal into a series of type \eqref{eq_series}, meaning that the coefficient sequence $(c_\lambda)_{\lambda \in \L}$ may not be square-summable \cite{LyuSei99}. This is a manifestation of the Balian-Low theorem. In order to obtain stable expansions of type \eqref{eq_series}, it is necessary that $G(g,\L)$ forms a (Gabor) frame for $\Lt[]$.

A Gabor system $\G(g,\L)$ is a frame for $\Lt[]$ if and only if there exist positive constants $0 < A \leq B < \infty$, called frame bounds, such that
\begin{equation}\label{eq_frame}
	A \norm{f}_2^2 \leq \sum_{\l \in \L} |\langle f, \pi(\l) g \rangle|^2 \leq B \norm{f}_2^2, \quad \forall f \in \Lt[].
\end{equation}
In case of the Gaussian (as well as for other ``suitably nice" functions), a Gabor system can only be a frame if it is redundant (overcomplete). This is achieved by increasing the (lower Beurling) density of the index set, which in case of a lattice $\L = M \Z^2$ is simply given by
\begin{equation}
	\delta(\L) = \frac{1}{|\det(M)|} \, .
\end{equation}

In case of a Gaussian window, the necessary density condition on the index set is already sufficient as proved in \cite{Lyu92}, \cite{Sei92}, \cite{SeiWal92}. In particular, a Gaussian Gabor system with lattice $\L$ is a frame if and only if $\delta(\L) > 1$. The result actually holds in a more general setting, namely for relatively separated point sets with lower Beurling density greater than 1.

There is a natural operator associated to a Gabor system, the frame operator;
\begin{equation}
	S_{g,\L} f = \sum_{\l \in \L} \langle f, \pi(\l) g \rangle \, \pi(\l) g.
\end{equation}
Its operator norm and the norm of its inverse are connected to the optimal frame bounds in \eqref{eq_frame} in the following way;
\begin{equation}
	\norm{S_{g,\L}}_{op} = B \qquad \textnormal{ and } \qquad \norm{S_{g,\L}^{-1}}_{op} = A^{-1}.
\end{equation}
The condition number of the frame operator is given by
\begin{equation}
	\cond(S_{g,\L}) = \frac{B}{A},
\end{equation}
where $A$ and $B$ are the sharp frame bounds, depending on the window $g$ and the lattice $\L$. Now, if the Gabor system constitutes a frame, the frame operator is invertible and the coefficients in \eqref{eq_series} can, e.g., be computed by using the canonical dual window $g^\circ = S_{g,\L}^{-1} g$;
\begin{equation}
	f = \sum_{\l \in \L} \langle f, \pi(\l) g^\circ \rangle \, \pi(\l) g.
\end{equation}

\section{The Conjecture of Strohmer and Beaver}
\begin{conjecture}[Strohmer and Beaver]\label{con_SB}
Consider the family of Gaussian Gabor systems $\G(g_0, \L)$ with fixed lattice density greater than 1, i.e., $\delta(\L) > 1$. Then, the condition number of the associated family of frame operators, $\cond(S_{g_0,\L}) = B/A$, is minimal, if and only if
\begin{equation}
	\L = \L_h = \delta^{-1/2} Q M_h \Z^2
\end{equation}
is a hexagonal lattice. Here, $\delta$ is the given density, $Q$ is an orthogonal matrix and
\begin{equation}
	M_h = \sqrt{\tfrac{2}{\sqrt{3}}}
	\left(
		\begin{array}{c c}
			1 & \frac{1}{2}\\
			0 & \frac{\sqrt{3}}{2}
		\end{array}
	\right).
\end{equation}
Considering the special case of rectangular (or separable) lattices of the form
\begin{equation}
	\L_{(\alpha, \beta)} = \alpha \Z \times \beta \Z =
	\left(
		\begin{array}{c c}
			\alpha & 0\\
			0 & \beta
		\end{array}
	\right) \Z^2,
\end{equation}
$(\alpha \beta)^{-1} > 1$ fixed, the square lattice ($\alpha = \beta$) minimizes the condition number of the frame operator.
\end{conjecture}

In shorter notation, the main claim in Conjecture \ref{con_SB} is that
\begin{equation}
	\cond(S_{g_0,\L_h}) \leq \cond(S_{g_0, \L}),
\end{equation}
with equality if and only if $\L$ is another (rotated) version of the hexagonal lattice.

\subsection{Heuristic Arguments and Proof by Intimidation}
For the heuristics, we first need to define the short-time Fourier transform of a function $f$ with respect to a window $g$;
\begin{equation}
	V_g f(x, \omega) = \int_{\R} f(t) \, \overline{g(t-x)} e^{-2 \pi i \omega t} \, dt = \langle f, \pi(\l) g \rangle,
\end{equation}
with $\l = (x, \omega) \in \R^2$. Now for $f = g = g_0$, we get
\begin{equation}
	V_{g_0}g_0(x, \omega) = e^{- \pi i x \omega} e^{-\tfrac{1}{2} \pi (x^2+\omega^2)}.
\end{equation}
The function $|V_g f|^2$ is called the spectrogram of $f$ with respect to the window $g$ and measures the time-frequency concentration of $f$ (with respect to $g$). In case of the Gaussian we get
\begin{equation}
	|V_{g_0}g_0(x,\omega)|^2 = e^{-\pi (x^2 + \omega^2)}.
\end{equation}
We see that this function is radial symmetric and most of its energy is concentrated in a disc. Since the optimal way to arrange discs in the plane is given by the hexagonal lattice, the first guess is that this is also the optimal way to arrange two-dimensional Gaussians (to be made precise below).

For $g \in \Lt[]$ and a lattice $\L \subset \R^2$, we set
\begin{equation}
	p_{g,\L}(z) = \sum_{\l \in \L}|V_g g (\l + z)|^2,
\end{equation}
which is $\L$-periodic in $z$. Then, for any window $g$ with $\norm{g}_2 = 1$, it follows from \eqref{eq_frame} that
\begin{equation}\label{eq_approx}
	A \leq \textnormal{ess} \inf_{z \in \R^2} p_{g,\L}(z) \quad \textnormal{ and } \quad \textnormal{ess} \sup_{z \in \R^2} p_{g,\L}(z) \leq B
\end{equation}
by considering not all $f \in \Lt[]$, but only all possible time-frequency shifted windows $\pi(z) g$, $z \in \R^2$. Now, the conjecture of Strohmer and Beaver on the smallest condition number is implied by the following, stronger conjecture.
\begin{conjecture}\label{con}
	For fixed density $\delta$, the lower frame bound of the Gaussian Gabor system $\G(g_0,\L)$ is uniquely maximized by the hexagonal lattice and the upper frame bound is uniquely minimized in this case.
	
Also, for fixed density $\delta$, for the separable (or rectangular) Gaussian Gabor system $\G(g_0, \alpha \Z \times \beta \Z)$, $(\alpha \beta)^{-1} = \delta$, the lower frame bound is uniquely maximized and the upper frame bound is uniquely minimized if and only if $\alpha = \beta = \delta^{-1/2}$.
\end{conjecture}
For special densities, the separable case in Conjecture \ref{con} was proven in 2017 \cite{FauSte17}. More recently, for special densities it was proven in \cite{Fau18} that the hexagonal lattice uniquely minimizes the upper frame bound. The only problem in Conjecture \ref{con} which is open for all densities is the problem of maximizing the lower frame bound among all lattices.

Furthermore, we note that for $\delta(\L) \in 2 \N$, the sharp frame bounds of the Gabor system $\G(g_0,\L)$ are given by (see \cite{Fau18note}, \cite{Jan96})
\begin{equation}\label{eq_sharp}
	A = \textnormal{ess} \inf_{z \in \R^2} \widetilde{p}_{g_0,\L^\circ}(z)
	\quad \textnormal{and} \quad
	B = \textnormal{ess} \sup_{z \in \R^2} \widetilde{p}_{g_0,\L^\circ}(z),
\end{equation}
where $\widetilde{p}_{g_0, \L^\circ}$ is the following Fourier series with Gaussian coefficients;
\begin{equation}
	\widetilde{p}_{g_0, \L^\circ}(z) = \delta \sum_{\l^\circ \in \L^\circ} e^{-\tfrac{\pi}{2} |\l^\circ|^2} e^{2 \pi i \sigma(\l^\circ, z)},
\end{equation}
where $\sigma(\l^\circ, z) = \l^\circ_1 z_2 - \l^\circ_2 z_1$ is the standard symplectic form, $\delta$ is the lattice density and $\L^\circ = \delta \L$ is the adjoint lattice. Concluding from \eqref{eq_sharp}, a first step towards a solution of Conjecture \ref{con} would be to see whether, for any density, the $\L$-periodic function $p_{g_0,\L^\circ}(z)$ assumes its largest minimum and its smallest maximum for the hexagonal lattice, which results in a problem for the heat kernel on a family of two-dimensional tori.

We note that, by the Poisson summation formula and the special choice of the window being $g_0$, we can connect $\widetilde{p}_{g_0, \L^\circ}$ to the function $p_{g_0,\L}$ and the optimality problem for 2-dimensional Gaussians. It is actually not hard to show, by using the triangle inequality, that for any lattice $\L^\circ$ we have
\begin{equation}
	\widetilde{p}_{g_0,\L^\circ}(z) \leq \widetilde{p}_{g_0,\L^\circ}(0).
\end{equation}
The result of Montgomery on minimal theta functions \cite{Mon88} states that for fixed lattice density $\delta$
\begin{equation}
	\widetilde{p}_{g_0, \L^\circ_h}(0) \leq \widetilde{p}_{g_0, \L^\circ}(0)
\end{equation}
with equality if and only if $\L^\circ$ is another hexagonal lattice. For $\delta \in 2 \N$ fixed, this is equivalent to the result that the upper frame bound is minimal if and only if the lattice is hexagonal \cite{Fau18}. Without going into the details, we note that the commutation relations \eqref{eq_commrel} are the reason why we only get the equivalence for even lattice densities.

If one could show that, for any $\delta$ the minimum of $\widetilde{p}_{g_0, \L^\circ}$ is maximal if and only if the lattice is hexagonal, Conjecture \ref{con} would be proved for even lattice densities, implying that the conjecture of Strohmer and Beaver is true for even lattice densities. A major issue is that locating the minimum is not as easy as locating the maximum and, furthermore, numerical investigations show that the location of the minimum also depends on $\delta$ (see also \cite{Bae98}).

We close this section with the promised ``proof" of the Strohmer and Beaver conjecture by intimidation.

(Claim): We claim that Conjecture \ref{con} is true and, hence, for fixed lattice density the condition number of a Gaussian Gabor frame operator is minimal only for a hexagonal lattice.

(``Proof"): Which other lattice should yield the minimal condition number? $\square$

\section{Extremal Geometries}
We start this section with the celebration of the $90^{th}$ birthday of a theorem and a related open problem by Landau \cite{Lan29}, stated in 1929.
\begin{theorem}[Landau, 1929]\label{thm_Landau}
	Let $f: \D \to \C$ be a holomorphic map from the open unit disc $\D$ to the complex plane $\C$ with the property $|f'(0)| = 1$. Then, there exists an absolute constant $\mathcal{L} > 0$ such that an open disc $D_\mathcal{L}$ of radius $\mathcal{L}$ is contained in the image of $f(\D)$.
\end{theorem}
Landau's problem is to find the exact value of the constant $\mathcal{L}$, which can be defined in the following way;
\begin{align}
	\ell(f) & = \sup \{r \in \R_+ \mid D_r \subset f(\D), f \textnormal{ as in Theorem \ref{thm_Landau}} \},\\
	\mathcal{L} & = \inf \{ \ell(f) \mid f \textnormal{ as in Theorem \ref{thm_Landau}} \}.
\end{align}
We note that the problem of finding the exact value of $\mathcal{L}$ is invariant under translation and rotation, just as the problems stated in Conjecture \ref{con}. We have the following estimates on $\mathcal{L}$;
\begin{equation}
	\frac{1}{2} < \mathcal{L} \leq \mathcal{L}_+ = \frac{\Gamma\left( \tfrac{1}{3} \right) \Gamma\left( \tfrac{5}{6} \right)}{\Gamma\left( \tfrac{1}{6} \right)} = 0.543259 \ldots \; .
\end{equation}
The value for $\mathcal{L}_+$ was established in 1943 by Rademacher \cite{Rad43}, who also mentioned that the same value was already derived by Robinson in 1937, but this work was not published. The value $\mathcal{L}_+$ was derived by constructing (and properly scaling) the universal covering map $\phi$ of a once-punctured hexagonal torus. The once-punctured hexagonal torus can be identified with the complex plane minus a hexagonal lattice $\L_h$;
\begin{equation}
	\mathbb{T}^2_h \cong \C \big\backslash \L_h.
\end{equation}
The map $\phi$ is constructed as follows.
\begin{figure}[ht]
	\centering{
	\begin{minipage}{.45\textwidth}
		\subfigure[Tessellation of the unit disc with hyperbolic triangles.]{\includegraphics[width=.4\textwidth]{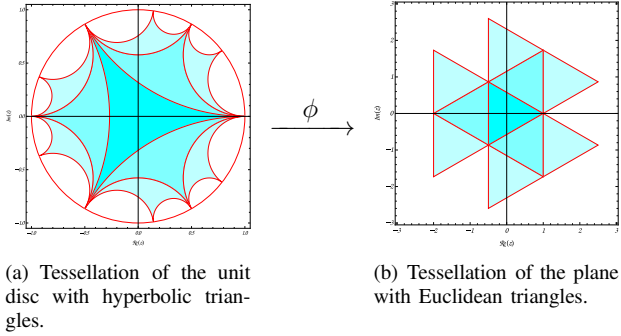}}
		\hfill \raisebox{1.5cm}{$\xrightarrow{\quad \displaystyle{\phi} \quad}$} \hfill
		\subfigure[Tessellation of the plane with Euclidean triangles.]{\includegraphics[width=.4\textwidth]{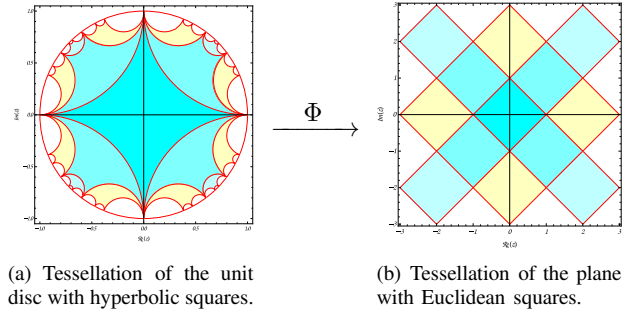}}
	\end{minipage}
	}
	\caption{Constructing the map $\phi$.}\label{fig_hex}
\end{figure}

One starts with a map $\phi_0$, mapping the unit disc to a hyperbolic equilateral triangle. As a second step, in the same manner one constructs a map $\phi_{1/3}$ from the unit disc to a Euclidean equilateral triangle. By composing the inverse map $\phi_0^{-1}$ with the map $\phi_{1/3}$ one maps the hyperbolic triangle to the Euclidean triangle. Finally, the map $\phi$ is constructed by using successive reflections in the unit disc (with the Poincare metric) and the plane, yielding a universal covering map of $\C \big \backslash \L_h$. The points in $\L_h$ become branching points of infinite order and, hence, $\phi$ is not holomorphic at these points. The largest disc that can be placed in $\C \big \backslash \L_h$ is the circumcircle of the constructed Euclidean triangle. The process of constructing this universal covering map $\phi$ is illustrated in Figure \ref{fig_hex}.

\begin{figure}[ht]
	\centering{
	\begin{minipage}{.45\textwidth}
		\subfigure[Tessellation of the unit disc with hyperbolic squares.]{\includegraphics[width=.4\textwidth]{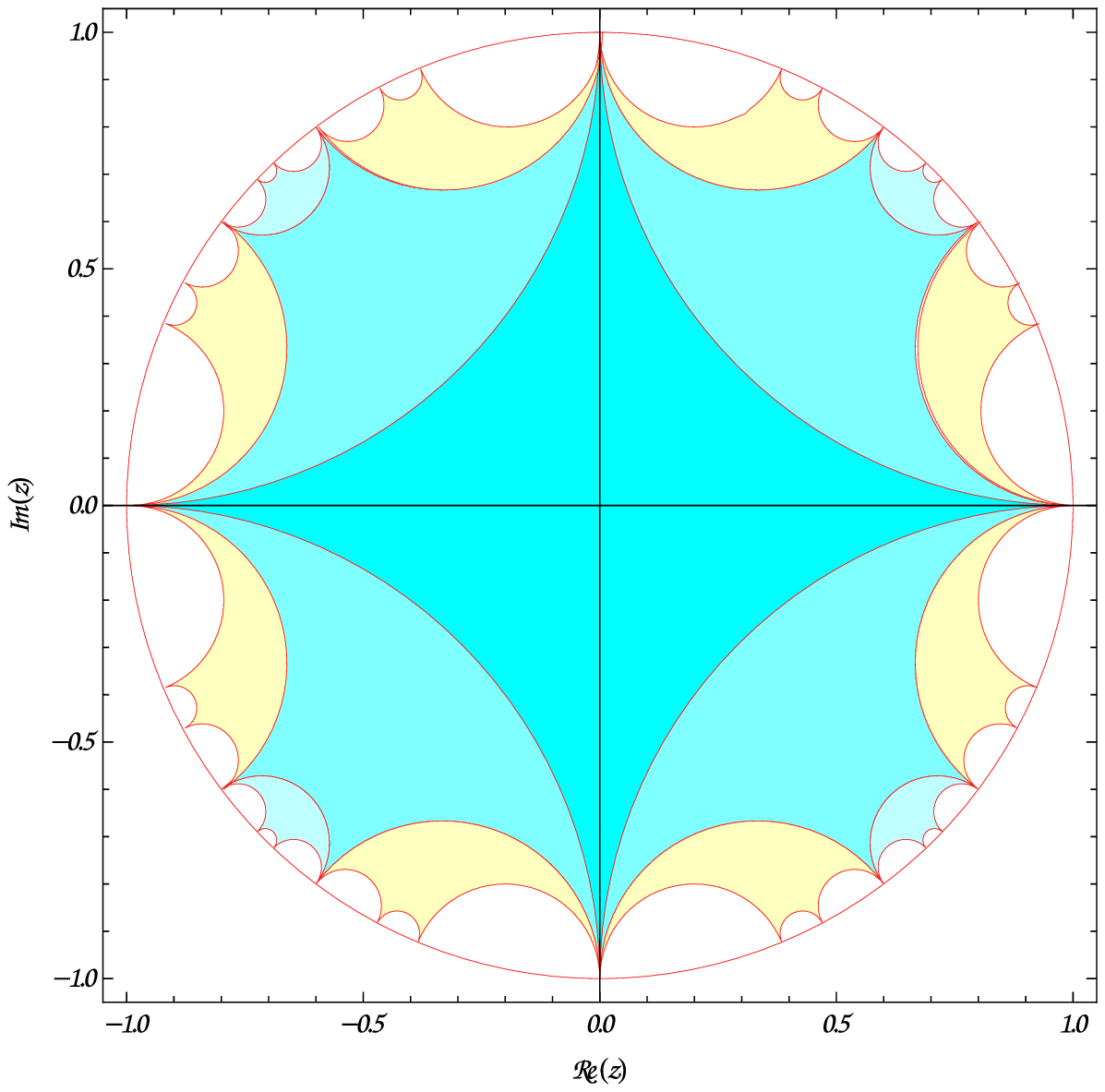}}
		\hfill \raisebox{1.5cm}{$\xrightarrow{\quad \displaystyle{\Phi} \quad}$} \hfill
		\subfigure[Tessellation of the plane with Euclidean squares.]{\includegraphics[width=.4\textwidth]{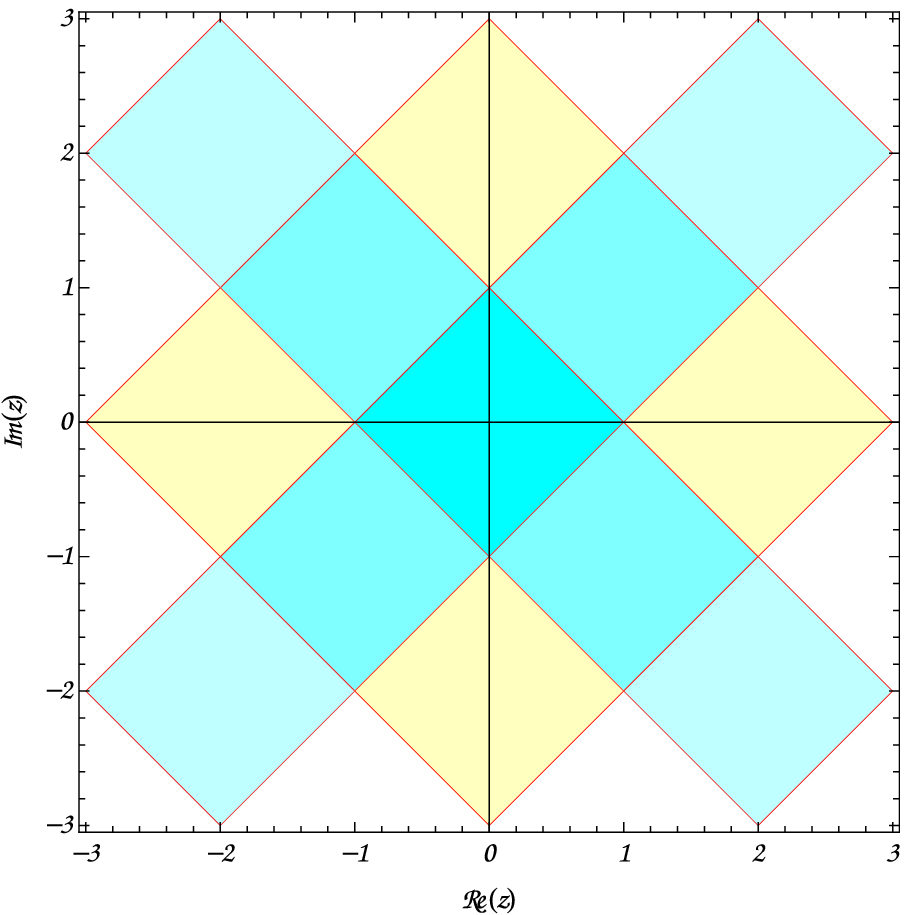}}
	\end{minipage}
	}
	\caption{Constructing the map $\Phi$.}\label{fig_square}
\end{figure}

By only considering universal covering maps of rectangular tori, the expected solution to the rectangular Landau problem is, of course, given by (scaling) the universal covering map $\Phi$ of the once-punctured square torus $\mathbb{T}^2 \cong \C \big\backslash \Z^2$, illustrated in Figure \ref{fig_square}. Note that the described mappings are not one-to-one.

The rectangular problem was investigated in \cite{Bae05} and \cite{Ere11}. The conjectured exact value of the rectangular Landau constant is
\begin{equation}\label{eq_Ls}
	\mathcal{L}_\square = \frac{\Gamma\left( \tfrac{1}{2} \right) \Gamma\left( \tfrac{3}{4} \right)}{\Gamma\left( \tfrac{1}{4} \right)} = 0.59907 \ldots \; .
\end{equation}

The sharp frame bounds of Gabor systems with the standard Gaussian window $g_0(t) = 2^{1/4} e^{-\pi t^2}$ and the hexagonal and square lattice of density 2 have been computed in \cite{Fau18}. The exact values of the lower frame bound of the Gaussian Gabor frame for the hexagonal and square lattice are
\begin{align}
	\mathcal{A} & = 2 \sum_{k,l \in \Z} e^{-\pi \tfrac{2}{\sqrt{3}} \left(k^2 + k l + l^2\right)} e^{2 \pi i \left(\tfrac{k}{3}-\tfrac{l}{3} \right)} = 1.84074 \ldots \; ,\\
	\mathbb{A} & = 2 \sum_{k,l \in \Z} e^{-\pi \left( k^2 + l^2 \right)} e^{2 \pi i \left( \tfrac{k}{2} - \tfrac{l}{2} \right)} = 1.66925 \ldots \; ,
\end{align}
respectively. For more details on how to compute sharp frame bounds for certain (integer) densities, we refer to \cite{Fau18note} and \cite{Jan96}.

By using results going back to Ramanujan and Gauss, it is shown in \cite{Fau19} that
\begin{align}\label{eq_constants}
	\mathbb{A} = \mathcal{L}_\square^{-1} \qquad \textnormal{ and } \qquad \mathcal{A} = \mathcal{L}_+^{-1} \; .
\end{align}
We will sketch the proof of \eqref{eq_constants} for $\mathbb{A}$, which uses Gauss' hypergeometric function ${}_2F_1$ and its connection to the Gamma function, established by Gauss, as well as the connection to theta functions, established by Ramanujan. For details we refer to the textbook of Berndt \cite[Chap.~17]{BerIII}.

First, note that $\mathbb{A}$ is expressible by means of Jacobi's theta functions;
\begin{equation}\label{eq_As}
	\mathbb{A} = 2 \, \theta_4(i)^2, \qquad \textnormal{ where } \qquad
	\theta_4(\tau) = \sum_{k_\in \Z} (-1)^k e^{\pi i \tau k^2}, 
\end{equation}
with $\tau \in \mathbb{H}$, the upper half plane. Also, for the square lattice of density 2 the upper frame bound, we denote it by $\mathbb{B}$, is expressible in a similar manner in terms of theta functions;
\begin{equation}
	\mathbb{B} = 2 \, \theta_3(i)^2, \qquad \theta_3(\tau) = \sum_{k_\in \Z} e^{\pi i \tau k^2}, \, \tau \in \mathbb{H}.
\end{equation}
Furthermore, we have
\begin{equation}\label{eq_modulus}
	\frac{\theta_4(i)^4}{\theta_3(i)^4} = \frac{1}{2},
\end{equation}
which is a well-known result (note that the value in \eqref{eq_modulus} yields the inverse of the squared condition number). Next, we introduce Gauss' hypergeometric function $_2{}F_1$;
\begin{equation}
	_2{}F_1(a,b;c;z) = \sum_{k = 0}^\infty \frac{(a)_k \, (b)_k}{(c)_k} \frac{z^k}{k!}, \qquad \textnormal{for }|z| < 1.
\end{equation}
For a number $w \in \C$ and $k \in \Z$, $(w)_k$ denotes the rising factorial, which is given as the ratio of Gamma functions;
\begin{equation}
	(w)_k = \frac{\Gamma(w+k)}{\Gamma(w)}.
\end{equation}
Gauss established the result
\begin{equation}\label{eq_Gauss}
	_2{}F_1\left(a,b;\tfrac{1}{2}(1+a+b);\tfrac{1}{2} \right) = \frac{\Gamma(\tfrac{1}{2}) \Gamma(\tfrac{1}{2}(1+a+b))}{\Gamma(\tfrac{1}{2}(1+a)) \Gamma(\tfrac{1}{2}(1+b))},
\end{equation}
and Ramanujan found the connection
\begin{equation}\label{eq_Ramanujan}
	_2{}F_1\left( \tfrac{1}{2}, \tfrac{1}{2}; 1; 1-\tfrac{\theta_4(\tau)^4}{\theta_3(\tau)^4} \right) = \theta_3(\tau)^2.
\end{equation}
Combining equations \eqref{eq_Ls}, \eqref{eq_As}, \eqref{eq_modulus}, \eqref{eq_Gauss} and \eqref{eq_Ramanujan} finally leads to the first equality in \eqref{eq_constants}. For further reading and more details we refer to \cite{Fau19} and the references therein.

The result for $\mathcal{A}$ follows in a similar manner by using Ramanujan's ``corresponding theory" \cite[Chap.~33]{BerV} and cubic analogues of Jacobi's theta functions \cite{BorBor91}.

Finally, we note that the problem of finding the exact value of Landau's constant and maximizing the lower frame bound heuristically have a lot in common. Both problems are invariant under translation and rotation and in both cases we fix a characteristic number. For Landau's problem it is the modulus of the derivative at the origin and in Conjecture \ref{con} it is the density of the lattice. Also, Landau's problem can be reduced to a problem for discrete subsets of $\C$ \cite{Bae98}. So, if $f$ is a universal covering map of $\C \backslash \L$ and $\L$ allows place for a large disc, so-to-say has a large hole, then the Gabor system $\G(g_0,\L)$ cannot have a large lower frame bound as one can find a function (e.g.~a time-frequency shifted Gaussian) essentially concentrated in this hole. Therefore, for this particular function the middle expression in  the frame inequality \eqref{eq_frame} will be small, forcing the lower frame bound to be small.

The given heuristic arguments together with \eqref{eq_constants} suggest that the conjecture of Strohmer and Beaver seems to be a rather deep mathematical problem. Also, we note that in Theorem \ref{thm_Landau} we can easily replace the assumption $|f'(0)| = 1$ by $|f'(0)| = K$, $K > 0$, as the scaling constant $K$ enters the problem linearly. This means that the problem is invariant under scaling and a solution for $K = 1$ already gives a solution for any $K \in \R_+$. On the other hand, if we solved the lower frame bound problem in Conjecture \ref{con} for density $\delta = 2$, it is not clear that we already solved the problem for any $\delta > 1$. From this point of view, solving Conjecture \ref{con} in full generality might be even harder than solving the Landau problem.

\section*{Acknowledgments}
\footnotesize{The author was supported by an Erwin--Schrödinger Fellowship of the Austrian Science Fund (FWF): J4100-N32. The presented results have partially been derived during the first stage of the fellowship, which the author spent with the Analysis Group at the Department of Mathematical Sciences, NTNU Trondheim, Norway.}

\end{document}